\documentclass[preprint,11pt,times]{elsarticle}

\usepackage{setspace}
\usepackage{geometry}
\usepackage{url}
\usepackage{graphicx,graphics}
\usepackage{subfigure}
\usepackage{epsfig}
\usepackage{psfrag}
\usepackage{hyperref}
\usepackage{amsmath,amssymb,amsfonts,amsthm}
\usepackage{mathrsfs}
\usepackage{color}
\usepackage{float}
\usepackage{natbib}
\usepackage{tabularx}
\usepackage{multirow}

\theoremstyle{remark}
\newtheorem{thm}{Theorem}[section]
\newtheorem{rmk}[thm]{Remark}

\newcommand{\Eref}[1]{Equation (\ref{#1})}
\newcommand{\fref}[1]{Figure (\ref{#1})}
\newcommand{\Erefs}[1]{Equations (\ref{#1})}

\newcommand{\xx}{\mathbf{x}}
\newcommand{\KK}{\mathbf{K}}
\newcommand{\bm}{\mathbf{M}}
\newcommand{\bn}{\mathbf{N}}
\newcommand{\DD}{\mathbf{D_b}}

\newcommand{\bveps}{\boldsymbol{\varepsilon}}

\newcommand{\BB}{\mathbf{B}}

\newcommand{\rmd}{\mathrm{d}}

\begin{document}

\begin{frontmatter}

\title{A parametric study on the buckling of functionally graded material plates with internal discontinuities using the partition of unity method}       

\author[unsw]{S Natarajan \corref{cor1}\fnref{fn1}}
\author[iisc]{S Chakraborty}
\author[india]{M Ganapathi}
\author[india2]{M Subramanian}

\cortext[cor1]{Corresponding author}

\address[unsw]{School of Civil \& Environmental Engineering, The University of New South Wales, Sydney, Australia}
\address[iisc]{Department of Aerospace Engineering, Indian Institute of Science, Bangalore, India}
\address[india]{Stress \& DTA, IES-Aerospace, Mahindra Satyam Computer Services Ltd., Bangalore, India}
\address[india2]{Professor, Department of Automobile Engineering, Dr. Mahalingam College of Engineering and Technology, Coimbatore, India}

\fntext[fn1]{\url sundararajan.natarajan@gmail.com}

\begin{abstract}

In this paper, the effect of local defects, viz., cracks and cutouts on the buckling behaviour of functionally graded material plates subjected to mechanical and thermal load is numerically studied. The internal discontinuities, viz., cracks and cutouts are represented independent of the mesh within the framework of the extended finite element method and an enriched shear flexible 4-noded quadrilateral element is used for the spatial discretization. The properties are assumed to vary only in the thickness direction and the effective properties are estimated using the Mori-Tanaka homogenization scheme. The plate kinematics is based on the first order shear deformation theory. The influence of various parameters, viz., the crack length and its location, the cutout radius and its position, the plate aspect ratio and the plate thickness on the critical buckling load is studied. The effect of various boundary conditions is also studied. The numerical results obtained reveal that the critical buckling load decreases with increase in the crack length, the cutout radius and the material gradient index. This is attributed to the degradation in the stiffness either due to the presence of local defects or due to the change in the material composition.
	
\end{abstract}


\end{frontmatter}


\section{Introduction}

The functionally graded materials (FGMs) are new class of engineered materials characterized by \emph{smooth and continuous transition} in properties from one surface to another~\cite{koizumi1993}. As a result, the FGMs are macroscopically homogeneous and are preferred over the laminated composites for structural integrity. The tunable thermo-mechanical property of the FGM has attracted researchers to study the static and the dynamic behaviour of structures made of FGM under mechanical~\cite{Zenkour2006,Zenkour2007,reddy2000,singhat.prakash2011} and thermal loading~\cite{natarajanbaiz2011, praveenreddy1998, dailiu2005, ganapathiprakash2006a, janghorbanazare2011, zenkourmashat2010, zhaolee2009}. Praveen \textit{et al.,}~\cite{praveenreddy1998} and Reddy \textit{et al.,}~\cite{reddychin2007} studied the thermo-elastic response of ceramic-metal plates using first order shear deformation theory (FSDT) coupled with 3D heat conduction equation. Their study concluded that the structures made up of FGM with ceramic rich side exposed to elevated temperatures are susceptible to buckling due to the through thickness temperature variation. The buckling of skewed FGM plates under mechanical and thermal loads were studied in ~\cite{ganapathiprakash2006a,ganapathiprakash2006} employing the FSDT and by using the shear flexible quadrilateral element. Efforts has also been made to study the mechanical behaviour of FGM plates with geometrical imperfection \cite{shariateslami2006}. Saji et. al. \cite{sajivarughese2008} studied thermal buckling of FGM plates with material properties dependent on both the composition and temperature. They found that the critical buckling temperature decreases when material properties are considered to be a function of temperature. Ganapati \textit{et al.,}~\cite{ganapathiprakash2006a} studied the buckling of FGM skewed plate under thermal loading. Efforts has also been made to study the mechanical behaviour of FGM plates with geometrical imperfection \cite{shariateslami2006}. FGM plates or in general plate structures, may develop flaws during manufacturing or after they have been put into service. Hence it is important to understand the response of a FGM plate with an internal flaw. It is known that cracks or local defects affect the response of a structural member. This is because, the presence of the crack introduces local flexibility and anisotropy. The vibration of cracked FGM structures are fairly dealt in the literature. Kitipornchai \textit{et al.,}~\cite{kitipornchaike2009} studied nonlinear vibration of edge cracked functionally graded Timoshenko beams using Ritz method. Yang \textit{et al.,}~\cite{yanghao2010} studied the nonlinear dynamic response of a functionally graded plate with a through-width crack based on Reddy's third-order shear deformation theory. Dolbow and Gosz~\cite{dolbowgosz2002} employed the extended finite element method (XFEM) to compute mixed mode stress intensity factors for a crack in a functionally graded material. Natarajan \textit{et al.,}~\cite{natarajanbaiz2011a,natarajan2011} studied the influence of cracks on the vibration and mechanical buckling of functionally graded material plates. The above list is no way comprehensive and interested readers are referred to the literature and references therein and a recent review paper by Jha and Kant~\cite{jhakant2013} on FGM plates. To the author's knowledge, the influence of the presence of an internal flaw, viz., cracks and cutouts on the critical buckling load and critical buckling temperature has not been studied earlier.

In this paper, we study the buckling behaviour of FGM plates with local defects, viz., cracks and cutouts. In this study, cracks and cutouts are considered as internal flaw. A structured quadrilateral mesh is employed and the local defects are represented independent of the underlying finite element (FE) mesh by enriching the FE approximation basis with additional functions. An enriched shear flexible 4-noded element proposed in~\cite{natarajanbaiz2011a,natarajan2011} is used for this study. The influence of various geometric parameters, viz., the plate aspect ratio, the thickness of the plate, the crack length, the crack orientation and location, the cutout radius, the support conditions and the gradient index on the critical buckling load is numerically studied.

The paper is organized as follows, the next section will give a brief over of Reissner-Mindlin plate theory and an introduction to FGM. Section \ref{spatdisct} discusses the spatial discretization within XFEM framework and numerical integration over enriched elements. Section~\ref{numexample} presents results for the buckling analyses of FGM plates with geometric defects(cracks) and material discontinuity(cutouts), followed by concluding remarks in the last section.

\section{Theoretical Formulation} \label{theory}
\subsection{Functionally graded material}

\begin{figure}
\centering
\input{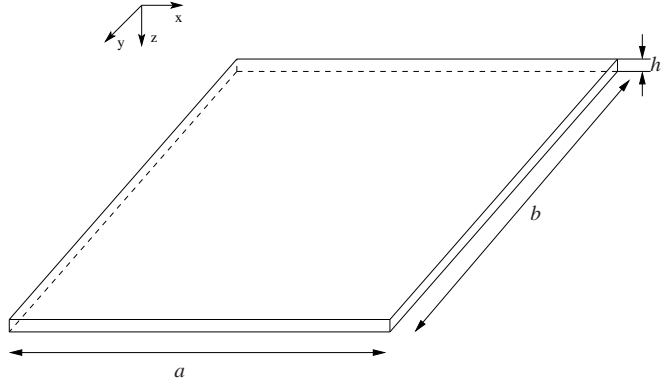}
\caption{Coordinate system of a rectangular FGM plate}
\label{fig:platefig}
\end{figure}

A functionally graded material (FGM) rectangular plate (length $a$, width $b$ and thickness $h$), made by mixing two distinct material phases: a metal and ceramic is considered with coordinates $x,y$ along the in-plane directions and $z$ along the thickness direction (see \fref{fig:platefig}). The material on the top surface $(z=h/2)$ of the plate is ceramic and is graded to metal at the bottom surface of the plate $(z=-h/2)$ by a power law distribution. The homogenized material properties are computed using the Mori-Tanaka Scheme~\cite{Benvensite1987,reddy2000,Qian2004a}. 

\subsubsection*{Estimation of mechanical and thermal properties}
\begin{figure}
\includegraphics{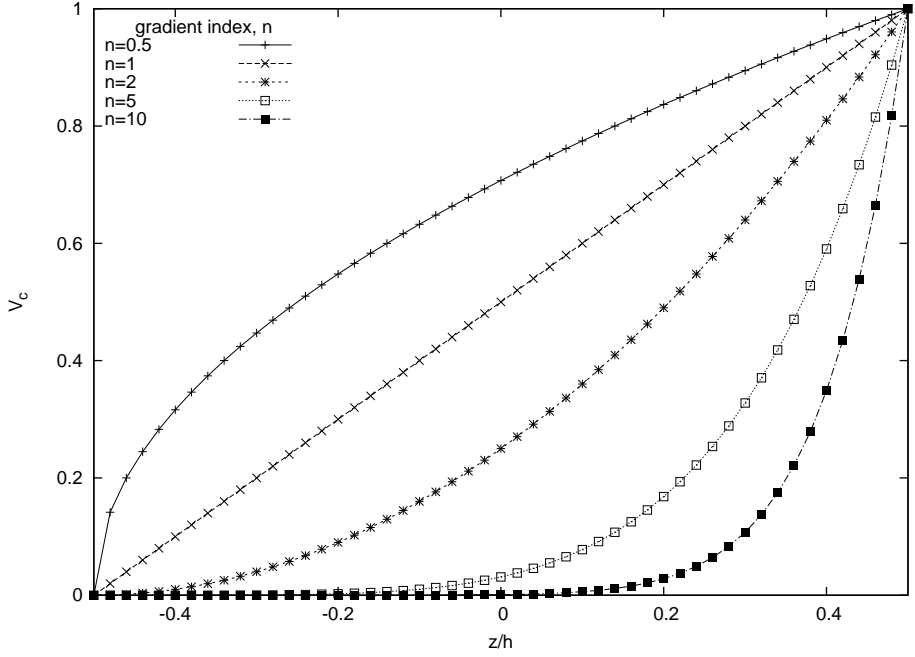}
\caption{Through thickness variation of volume fraction}
\label{fig:volfrac}
\end{figure}

Based on the Mori-Tanaka homogenization method, the effective bulk modulus $K$ and shear modulus $G$ of the FGM are evaluated as~\cite{Benvensite1987,Qian2004a}
\begin{eqnarray}
{K - K_m \over K_c - K_m} &=& {V_c \over 1+(1-V_c){3(K_c - K_m) \over 3K_m + 4G_m}} \nonumber \\
{G - G_m \over G_c - G_m} &=& {V_c \over 1+(1-V_c){(G_c - G_m) \over G_m + f_1}}
\label{eqn:bulkshearmodulus}
\end{eqnarray}
where
\begin{equation}
f_1 = {G_m (9K_m + 8G_m) \over 6(K_m + 2G_m)}
\end{equation}
Here, $V_i~(i=c,m)$ is the volume fraction of the phase material. The subscripts $c$ and $m$ refer to the ceramic and metal phases, respectively. The volume fractions of the ceramic and metal phases are related by $V_c + V_m = 1$, and $V_c$ is expressed as
\begin{equation}
V_c(z) = \left( {2z + h \over 2h} \right)^n, \hspace{0.2cm}  n \ge 0
\label{eqn:volFrac}
\end{equation}
where $n$ in \Eref{eqn:volFrac} is the volume fraction exponent, also referred to as the material gradient index. \fref{fig:volfrac} shows the variation of the volume fractions of ceramic and metal, respectively, in the thickness direction $z$ for the FGM plate. The top surface is ceramic rich and the bottom surface is metal rich. The effective Young's modulus $E$ and Poisson's ratio $\nu$ can be computed from the following expressions:
\begin{equation}
E = {9KG \over 3K+G} \hspace{0.3cm} \nu = {3K - 2G \over 2(3K+G)}
\label{eqn:young}
\end{equation}

The effective mass density $\rho$ is given by the rule of mixtures as $\rho = \rho_c V_c + \rho_m V_m$. The effective heat conductivity $\kappa_{\rm eff}$ and the coefficient of thermal expansion $\alpha_{\rm eff}$ is given by:
\begin{eqnarray}
\frac{\kappa_{\rm eff} - \kappa_m}{\kappa_c - \kappa_m} = \frac{V_c}{1 + V_m \frac{(\kappa_c - \kappa_m)}{3\kappa_m}} \nonumber \\
\frac{\alpha_{\rm eff} - \alpha_m}{\alpha_c - \alpha_m} = \frac{ \left( \frac{1}{K_{\rm eff}} - \frac{1}{K_m} \right)}{\left(\frac{1}{K_c} - \frac{1}{K_m} \right)}
\label{eqn:effAlpha}
\end{eqnarray}

\subsection*{Temperature distribution through the thickness}
The temperature variation is assumed to occur in the thickness direction only and the temperature field is considered to be constant in the $xy$-plane. In such a case, the temperature distribution along the thickness can be obtained by solving a steady state heat transfer equation

\begin{equation}
-{d \over dz} \left[ \kappa(z) {dT \over dz} \right] = 0, \hspace{0.5cm} T = T_c ~\textup{at}~ z = h/2;~~ T = T_m ~\textup{at} ~z = -h/2
\label{eqn:heat}
\end{equation}

The solution of \Eref{eqn:heat} is obtained by means of a polynomial series~\cite{Wu2004} as

\begin{equation}
T(z) = T_m + (T_c - T_m) \eta(z,h)
\label{eqn:tempsolu}
\end{equation}
where,
\begin{equation*}
\begin{split}
\eta(z,h) = {1 \over C} \left[ \left( {2z + h \over 2h} \right) - {\kappa_{cm} \over (n+1)\kappa_m} \left({2z + h \over 2h} \right)^{n+1} +  {\kappa_{cm} ^2 \over (2n+1)\kappa_m ^2 } \left({2z + h \over 2h} \right)^{2n+1}
-{\kappa_{cm} ^3 \over (3n+1)\kappa_m ^3 } \left({2z + h \over 2h} \right)^{3n+1} \right. \\ + 
\left. {\kappa_{cm} ^4 \over (4n+1)\kappa_m^4 } \left({2z + h \over 2h} \right)^{4n+1} 
- {\kappa_{cm} ^5 \over (5n+1)\kappa_m ^5 } \left({2z + h \over 2h} \right)^{5n+1} \right] ;
\end{split}
\label{eqn:heatconducres}
\end{equation*}

\begin{equation*}
\begin{split}
C = 1 - {\kappa_{cm} \over (n+1)\kappa_m} + {\kappa_{cm} ^2 \over (2 n+1)\kappa_m ^2} 
- {\kappa_{cm} ^3 \over (3n+1)\kappa_m ^3} \\ + {\kappa_{cm} ^4 \over (4n+1)\kappa_m ^4}
- {\kappa_{cm} ^5\over (5n+1)\kappa_m ^5}.
\end{split}
\end{equation*}

\subsection{Reissner-Mindlin Plate}
Using the Mindlin formulation, the displacements $u,v,w$ at a point $(x,y,z)$ in the plate (see \fref{fig:platefig}) from the medium surface are expressed as functions of the mid-plane displacements $u_o,v_o,w_o$ and independent rotations $\beta_x,\beta_y$ of the normal in $yz$ and $xz$ planes, respectively, as
\begin{eqnarray}
u(x,y,z,t) &=& u_o(x,y,t) + z \beta_x(x,y,t) \nonumber \\
v(x,y,z,t) &=& v_o(x,y,t) + z \beta_y(x,y,t) \nonumber \\
w(x,y,z,t) &=& w_o(x,y,t) 
\label{eqn:displacements}
\end{eqnarray}
where $t$ is the time. The strains in terms of mid-plane deformation can be written as:
\begin{equation}
\bveps  = \left\{ \begin{array}{c} \bveps_p \\ \mathbf{0} \end{array} \right \}  + \left\{ \begin{array}{c} z \bveps_b \\ \bveps_s \end{array} \right\} 
\label{eqn:strain1}
\end{equation}
The midplane strains $\bveps_p$, the bending strains $\bveps_b$ and the shear strains $\varepsilon_s$ in \Eref{eqn:strain1} are written as:
\begin{equation}
\renewcommand{\arraystretch}{1}
\bveps_p = \left\{ \begin{array}{c} u_{o,x} \\ v_{o,y} \\ u_{o,y}+v_{o,x} \end{array} \right\}, \hspace{0.2cm}
\renewcommand{\arraystretch}{1}
\bveps_b = \left\{ \begin{array}{c} \beta_{x,x} \\ \beta_{y,y} \\ \beta_{x,y}+\beta_{y,x} \end{array} \right\}, \hspace{0.2cm}
\renewcommand{\arraystretch}{1}
\bveps_s = \left\{ \begin{array}{c} \beta _x + w_{o,x} \\ \beta _y + w_{o,y} \end{array} \right\}.
\label{eqn:platestrain}
\end{equation}
where the subscript `comma' represents the partial derivative with respect to the spatial coordinate succeeding it. The membrane stress resultants $\bn$ and the bending stress resultants $\bm$ can be related to the membrane strains, $\bveps_p$ and the bending strains $\bveps_b$ through the following constitutive relations:
\begin{eqnarray}
\bn &=& \left\{ \begin{array}{c} N_{xx} \\ N_{yy} \\ N_{xy} \end{array} \right\} = \mathbf{A} \bveps_p + \BB \bveps_b - \bn^\textup{th} \nonumber \\
\bm &=& \left\{ \begin{array}{c} M_{xx} \\ M_{yy} \\ M_{xy} \end{array} \right\} = \BB \bveps_p + \DD \bveps_b - \bm^\textup{th} 
\end{eqnarray}
where the matrices $\mathbf{A} = A_{ij}, \BB= B_{ij}$ and $\DD = D_{ij}; (i,j=1,2,6)$ are the extensional, the bending-extensional coupling and the bending stiffness coefficients and are defined as
\begin{equation}
\left\{ A_{ij}, ~B_{ij}, ~ D_{ij} \right\} = \int_{-h/2}^{h/2} \overline{Q}_{ij} \left\{1,~z,~z^2 \right\}~dz
\end{equation}
The thermal stress resultant, $\bn^{\textup{th}}$ and the moment resultants $\bm^{\textup{th}}$ are:
\begin{eqnarray}
\bn^{\textup{th}} &=& \left\{ \begin{array}{c} N_{xx}^{\textup{th}} \\ N_{yy}^{\textup{th}} \\ N_{xy}^{\textup{th}} \end{array} \right\} = \int_{-h/2}^{h/2} \overline{Q}_{ij} \alpha(z,T) \left\{ \begin{array}{c} 1 \\ 1 \\ 0 \end{array} \right\} ~ \Delta T(z) ~ dz \nonumber \\
\bm^{\textup{th}} &=& \left\{ \begin{array}{c} M_{xx}^{\textup{th}} \\ M_{yy}^{\textup{th}} \\ M_{xy}^{\textup{th}} \end{array} \right\} = \int_{-h/2}^{h/2} \overline{Q}_{ij} \alpha(z,T) \left\{ \begin{array}{c} 1 \\ 1 \\ 0 \end{array} \right\} ~z~\Delta T(z) ~ dz
\end{eqnarray}
where the thermal coefficient of expansion $\alpha(z,T)$ is given by \Eref{eqn:effAlpha} and $\Delta T(z) = T(z)-T_o$ is the temperature rise from the reference temperature $T_o$ at which there are no thermal strains. Similarly, the transverse shear force $\mathbf{Q} = \{Q_{xz},Q_{yz}\}$ is related to the transverse shear strains $\varepsilon_s$ through the following equation
\begin{equation}
\mathbf{Q} = \mathbf{E} \bveps_s
\end{equation}
where $\mathbf{E} = E_{ij} = \int_{-h/2}^{h/2} \overline{Q}_{ij} \upsilon_i \upsilon_j~dz;~ (i,j=4,5)$ are the transverse shear stiffness coefficients and $\upsilon_i, \upsilon_j$ are the transverse shear coefficients for non-uniform shear strain distribution through the plate thickness. The stiffness coefficients $\overline{Q}_{ij}$ are defined as:
\begin{eqnarray}
\overline{Q}_{11} = \overline{Q}_{22} = {E(z,T) \over 1-\nu^2}; \hspace{0.2cm} \overline{Q}_{12} = {\nu E(z,T) \over 1-\nu^2}; \hspace{0.2cm} \overline{Q}_{16} = \overline{Q}_{26} = 0; \nonumber \\
\overline{Q}_{44} = \overline{Q}_{55} = \overline{Q}_{66} = {E(z,T) \over 2(1+\nu) }; \hspace{0.2cm} \overline{Q}_{45}=\overline{Q}_{54} = 0.
\end{eqnarray}
where the modulus of elasticity $E(z,T)$ and Poisson's ratio $\nu$ are given by \Eref{eqn:young}. The strain energy function $U$ is given by:
\begin{equation}
U(\boldsymbol{\delta}) = {1 \over 2} \int_{\Omega} \left\{ \bveps_p^{\textup{T}} \mathbf{A} \bveps_p + \bveps_p^{\textup{T}} \mathbf{B} \bveps_b + 
\bveps_b^{\textup{T}} \mathbf{B} \bveps_p + \bveps_b^{\textup{T}} \mathbf{D} \bveps_b +  \bveps_s^{\textup{T}} \mathbf{E} \bveps_s -  \bveps_p^{\textup{T}} \bn ^{\rm th}- \bveps_b^{\textup{T}} \bm^{\rm th} \right\}~ \rmd\Omega
\label{eqn:potential}
\end{equation}


where $\boldsymbol{\delta} = \{u,v,w,\beta_x,\beta_y\}$ is the vector of the degrees of freedom associated to the displacement field in a finite element discretization. The FGMs are most suited in high temperature environment, where the plate is subjected to a temperature field and this in turn results in in-plane stress resultants $(N_{xx}^{\textup{th}}, N_{yy}^{\textup{th}}, N_{xy}^{\textup{th}})$. The external work due to the in-plane stress resultants $(N_{xx}^{\textup{th}}, N_{yy}^{\textup{th}}, N_{xy}^{\textup{th}})$ developed in the plate under the thermal load is

\begin{equation}
\begin{split}
V(\boldsymbol{\delta}) = \int_{\Omega} \left\{ {1 \over 2} \left[ N_{xx}^{\textup{th}} w_{o,x}^2 + N_{yy}^{\textup{th}} w_{o,y}^2 +
2 N_{xy}^{\textup{th}} w_{o,x} w_{o,y}\right] \right. \\
\left. + {h^2 \over 24} \left[ N_{xx}^{\textup{th}} \left( \theta_{x,x}^2 + \theta_{y,x}^2 \right) + N_{yy}^{\textup{th}} \left( \theta_{x,y}^2 + \theta_{y,y}^2 \right) + 2N_{xy}^{\textup{th}} \left( \theta_{x,x} \theta_{x,y} + \theta_{y,x} \theta_{y,y} \right) \right] \right\} ~\rmd \Omega
\end{split}
\label{eqn:potthermal}
\end{equation}
Substituting \Eref{eqn:potential} - (\ref{eqn:potthermal}) in Lagrange's equations of motion and following the procedure given in~\cite{Rajasekaran1973}, the following discretized equation is obtained:

\paragraph*{Mechanical Buckling}
\begin{equation}
\left( \KK + \lambda \KK_G \right) \boldsymbol{\delta} = \mathbf{0}
\end{equation}

\paragraph*{Thermal Buckling}
\begin{equation}
\left( \KK + \Delta T \KK_G \right) \boldsymbol{\delta} = \mathbf{0}
\end{equation}

where $\Delta T(=T_c - T_m)$ is the critical temperature difference, $\lambda$ is the critical buckling load and $\KK$, $\KK_G$ are the linear stiffness and the geometric stiffness matrices, respectively. The critical buckling load and the critical temperature difference is computed using a standard eigenvalue algorithm.

\section{Spatial discretization}\label{spatdisct}

\subsection{Element description and shear locking}
The plate element employed here is a $\mathcal{C}^o$ continuous shear flexible field consistent element with five degrees of freedom $(u_o,v_o,w_o,\beta_x,\beta_y)$ at four nodes in a 4-noded quadrilateral (QUAD-4) element. The displacement field within the element is approximated by:
\begin{equation}
\{ u_o^e,v_o^e,w_o^e,\beta_x^e,\beta_y^e\} = \sum\limits_{J=1}^4 N_J \{u_{oJ}, v_{oJ}, w_{oJ},\beta_{xJ},\beta_{yJ} \},
\end{equation} 
where $u_{oJ}, v_{oJ}, w_{oJ},\beta_{xJ},\beta_{yJ}$ are the nodal variables and $N_J$ are the shape functions for the bi-linear QUAD-4 element, given by:
\begin{eqnarray}
N_1(\xi,\eta) = \frac{1}{4} (1-\xi)(1-\eta), \hspace{1cm}
N_2(\xi,\eta) = \frac{1}{4} (1+\xi)(1-\eta) \nonumber \\
N_3(\xi,\eta) = \frac{1}{4} (1+\xi)(1+\eta), \hspace{1cm}
N_4(\xi,\eta) = \frac{1}{4} (1-\xi)(1+\eta).
\label{eqn:Q4InterFun}
\end{eqnarray}
where $-1 \le \xi \le 1$ and $-1 \le \eta \le 1$. If the interpolation functions, given by \Eref{eqn:Q4InterFun} for a QUAD-4 are used directly to interpolate the five variables $(u_o,v_o,w_o,\beta_x,\beta_y)$ in deriving the shear strains and the membrane strains, the element will lock and show oscillations in the shear and the membrane stresses. The oscillations are due to the fact that the derivative functions of the out-of plate displacement, $w_o$ do not match that of the rotations ($\beta_x, \beta_y$) in the shear strain definition, given by \Eref{eqn:platestrain}. To alleviate the locking phenomenon, the terms corresponding to the derivative of the out-of plate displacement, $w_o$ must be consistent with the rotation terms, $\beta_x$ and $\beta_y$. The present formulation, when applied to thin plates, exhibits shear locking. In this study, field redistributed shape functions are used to alleviate the shear locking.~\cite{natarajanbaiz2011a,natarajan2011} The field consistency requires that the transverse shear strains and the membrane strains must be interpolated in a consistent manner. Thus, the $\beta_x$ and $\beta_y$ terms in the expressions for the shear strain $\bveps_s$ have to be consistent with the derivative of the field functions, $w_{o,x}$ and $w_{o,y}$. If the element has edges which are aligned with the coordinate system $(x,y)$, the section rotations $\beta_x, \beta_y$ in the shear strain are approximated by~\cite{somashekarprathap1987}:
\begin{equation}
\beta_x^e = \sum\limits_{J=1}^4 \tilde{N}_{1J} \beta_{xJ}, \hspace{0.5cm}
\beta_y^e  = \sum\limits_{J=1}^4 \tilde{N}_{2J}  \beta_{yJ}.
\end{equation}
where $\beta_{xJ}$ and $\beta_{yJ}$ are the nodal variables, $\tilde{N}_{1J}$ and $\tilde{N}_{2J}$ are the substitute shape functions, given by~\cite{somashekarprathap1987}:
\begin{eqnarray}
\tilde{N}_{1}(\eta) = \frac{1}{4} \left[ \begin{array}{cccc} 1-\eta & 1-\eta & 1+\eta &  1+\eta \end{array} \right] \nonumber \\
\tilde{N}_{2}(\xi) = \frac{1}{4} \left[ \begin{array}{cccc} 1-\xi  & 1 +\xi & 1+\xi  & 1-\xi \end{array} \right].
\label{eqn:fieldredistribute}
\end{eqnarray}
It can be seen that the field redistributed shape functions, given by \Eref{eqn:fieldredistribute} are consistent with the derivative of the shape functions given by \Eref{eqn:Q4InterFun} used to approximate the out-of plate displacement, $w_o$. Note that, no special integration rule is required for evaluating the shear terms. 

\subsection{Representation of discontinuity surface}
The finite element framework requires the underlying finite element mesh to conform to the discontinuity surface. The recent introduction of implicit boundary definition-based methods, viz., the extended/generalized FEM (XFEM/GFEM), alleviates the shortcomings associated with the meshing of the discontinuity surface. In this study, the partition of unity framework is employed to represent the discontinuity surface independent of the underlying mesh.

\paragraph*{Cracks}
\fref{fig:xfem_crack} shows a structured quadrilateral mesh with the crack and the cutout represented independent of the mesh. To represent the crack, two set of functions are used: (a) a Heaviside function to capture the jump across the crack face and (b) a set of functions that span the asymptotic displacement fields. The displacement approximation is divided into two parts: (a) standard part and (b) enriched part and the following enriched approximation for the plate displacements and the section rotations are used:

\begin{figure}[htpb]
\centering
\includegraphics[angle=0,width=0.4\columnwidth]{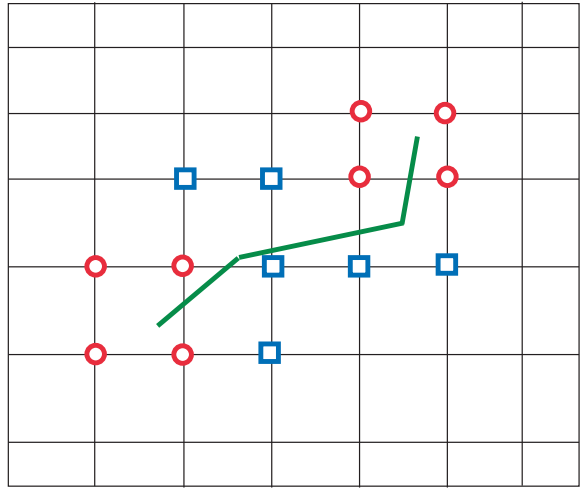}
\caption{A typical FE mesh with an arbitrary crack. `Squared' nodes are enriched with the heaviside function and `circled' nodes with the near tip functions, which allows representing cracks independent of the background mesh.}
\label{fig:xfem_crack}
\end{figure}

\begin{equation}
\begin{split}
(u^h,v^h,w^h)\left(\xx\right) = \underbrace{ \sum_{I \in \mathcal{N}^{\rm{fem}}} N_I(\xx) (u_I^s,v_I^s,w_I^s)}_{\rm FEM} + \underbrace{ \sum_{J \in
\mathcal{N}^{\rm{c}}} N_J(\xx) H(\xx) (b_J^u,b_J^v,b_J^w)}_{\rm Enriched~part} +  \\ \underbrace{ \sum_{K \in \mathcal{N}^{\rm{f}}}N_K(\xx) \left(
\sum_{l=1}^5 (c_{Kl}^u,c_{Kl}^v,c_{Kl}^w) G_l (r,\theta) \right)}_{\rm Enriched~part}
\end{split}
\label{eqn:platexfem1}
\end{equation}
The section rotations are approximated by:
\begin{eqnarray}
\beta_x^h \left(\xx\right) = \underbrace{ \sum_{I \in \mathcal{N}^{\rm{fem}}} \tilde{N}_{1I}(\xx) \beta_{x_I}^s}_{\rm FEM} +
\underbrace{ \sum_{J \in \mathcal{N}^{\rm{c}}} \tilde{N}_{1J}(\xx) H(\xx) b_J^{\beta_x} +
\sum_{K \in \mathcal{N}^{\rm{f}}} \tilde{N}_{1K}(\xx) \left( \sum_{l=1}^4 c_{Kl}^{\beta_x} F_l
(r,\theta) \right)}_{\rm Enriched~part}, \nonumber \\
\beta_y^h \left(\xx\right) = \underbrace{ \sum_{I \in \mathcal{N}^{\rm{fem}}} \tilde{N}_{2I}(\xx) \beta_{y_I}^s}_{\rm FEM} +
\underbrace{ \sum_{J \in \mathcal{N}^{\rm{c}}} \tilde{N}_{2J}(\xx) H(\xx) b_J^{\beta_y} + 
\sum_{K \in \mathcal{N}^{\rm{f}}} \tilde{N}_{2K}(\xx) \left( \sum_{l=1}^4 c_{Kl}^{\beta_y} F_l
(r,\theta) \right)}_{\rm Enriched~part}. 
\label{eqn:platexfem2}
\end{eqnarray}
where $\mathcal{N}^{\rm{fem}}$ is a set of all the nodes in the finite element mesh, $\mathcal{N}^{\rm{c}}$ is a set of nodes that are enriched with the Heaviside function and $\mathcal{N}^{\rm{f}}$ is a set of nodes that are enriched with near-tip asymptotic fields. In \Erefs{eqn:platexfem1} and (\ref{eqn:platexfem2}), $(u_I^s,v_I^s,w_I^s,\beta_{x_I}^s,\beta_{y_I}^s)$ are the nodal unknown vectors associated with the continuous part of the finite element solution, $b_J$ is the nodal enriched degree of freedom vector associated with the Heaviside (discontinuous) function, and $c_{Kl}$ is the nodal enriched degree of freedom vector associated with the elastic asymptotic near-tip functions. The asymptotic functions, $G_l$ and $F_l$ in \Erefs{eqn:platexfem1} and (\ref{eqn:platexfem2}) are given by
(\cite{dolbowmoes2000}):
\begin{eqnarray}
G_l(r,\theta) \equiv \left\{ \sqrt{r} \sin\left(\frac{\theta}{2}\right),\sqrt[3]{r}
\sin\left(\frac{\theta}{2}\right), \sqrt[3]{r}
\cos\left(\frac{\theta}{2}\right), \sqrt[3]{r}
\sin\left(\frac{3\theta}{2}\right), \sqrt[3]{r}
\cos\left(\frac{3\theta}{2}\right) \right\}, \nonumber \\
F_l(r,\theta) \equiv \sqrt{r} \left\{
\sin\left(\frac{\theta}{2}\right),
\cos\left(\frac{\theta}{2}\right),
\sin\left(\frac{\theta}{2}\right)\sin\left(\theta\right),
\cos\left(\frac{\theta}{2}\right)\sin\left(\theta\right)\right\}.
\label{eqn:asymp}
\end{eqnarray}
Here $(r,\theta)$ are the polar coordinates in the local coordinate system with the origin at the crack tip.


\begin{rmk} As we are interested in the global behaviour of the cracked FGM plate, we propose to use the same enrichment functions. The role of these enrichment functions is to aid in representing the discontinuous surface independent of the mesh.
\end{rmk}

\paragraph*{Cutouts} In this study, a level set approach is followed to model the cutouts. The geometric interface (for example, the boundary of the cutout) is represented by the zero level curve $\phi \equiv \phi(\xx,t) = 0$. The interface is located from the value of the level set information stored at the nodes. The standard FE shape functions can be used to interpolate $\phi$ at any point $\xx$ in the domain as:
\begin{equation}
\phi(\xx) = \sum\limits_I N_I(\xx) \phi_I
\end{equation}
where the summation is over all the nodes in the connectivity of the elements that contact $\xx$ and $\phi_I$ are the nodal values of the level set function. For circular cutout, the level set function is given by:
\begin{equation}
\phi_I = || \xx_I - \xx_c|| - r_c
\end{equation}
where $\xx_c$ and $r_c$ are the center and the radius of the cutout. For an elliptical cutout oriented at an angle $\theta$, measured from the $x-$ axis the level set function is given by:
\begin{equation}
\phi_I = \sqrt{a_1(x_I-x_c)^2-a_2(x_I-x_c)(y_I-y_c) + a_3(y_I-y_c)^2} - 1 
\end{equation}
where
\begin{equation}
a_1 = \left( {\cos \theta \over d} \right)^2, \hspace{0.15cm} a_2 = 2\cos\theta \sin\theta\left( {1 \over d^2} - {1 \over e^2} \right), \hspace{0.15cm} a_3 = \left( {\sin \theta \over d} \right)^2 + \left( {\cos \theta \over e}\right)^2.
\end{equation}
where $d$ and $e$ are the major and minor axes of the ellipse and $(x_c,y_c)$ is the center of the ellipse. 

\subsection{Numerical integration over enriched elements}
A consequence of adding custom tailored enrichment functions to the FE approximation basis, which are not necessarily smooth functions is that, special care has to be taken in numerically integrating over the elements that are intersected by the discontinuity surface. The standard Gau\ss~quadrature cannot be applied in elements enriched by discontinuous terms, because Gau\ss~quadrature implicitly assumes a polynomial approximation. One potential solution for the purpose of numerical integration is by partitioning the elements into subcells (to triangles for example) aligned to the discontinuous surface in which the integrands are continuous and differentiable~\cite{belytschkoblack1999}. The other techniques that can be employed are Schwarz Christoffel Mapping~\cite{natarajanbordas2009,natarajanmahapatra2010}, Generalized quadrature~\cite{mousavisukumar2011} and Smoothed eXtended FEM~\cite{bordasnatarajan2011}. In the present study, a triangular quadrature with sub-division is employed along with the integration rules described in Table \ref{table:subcellgausspt}. For the elements that are not enriched, a standard 2 $\times$ 2 Gaussian quadrature rule is used.

\begin{table}[htpb]
\caption{Integration rules for enriched and non-enriched elements in the presence of a crack}
\centering
\begin{tabular}{lr}
\hline
Element Type & Gau\ss ~points\\
\hline
Non-enriched element & 4  \\
Tip element & 13 per triangle \\
Tip blending element & 16 \\
Split element & 3 per triangle  \\
Split blending element & 4 \\
Split-Tip blending element & 4 per triangle \\
\hline
\end{tabular}
\label{table:subcellgausspt}
\end{table}

\section{Numerical Examples} \label{numexample}
In this section, we study the influence of local defects, viz., cracks and cutouts on the mechanical and thermal buckling behaviour of FGM plates. We consider both thin and thick plates with two different boundary conditions, viz., all edges simply supported (SSSS) and all edges clamped (CCCC). The boundary conditions for simply supported and clamped cases are :

\noindent \emph{Simply supported boundary condition}: \\
\begin{equation}
u_o = w_o = \theta_y = 0 \hspace{0.2cm} ~\textup{on} \hspace{0.2cm}  x=0,a; \hspace{0.2cm}
v_o = w_o = \theta_x = 0 \hspace{0.2cm} ~\textup{on} \hspace{0.2cm}  y=0,b
\end{equation}

\noindent \emph{Clamped boundary condition}: \\
\begin{equation}
u_o = w_o = \theta_y = v_o = \theta_x = 0 \hspace{0.2cm} ~\textup{on} ~ x=0,a \hspace{0.2cm} \& \hspace{0.2cm}  y=0,b
\end{equation}

The FGM plate considered here consists of Aluminum (Al) and Zirconium dioxide (ZrO$_2$). The material is considered to be temperature independent. The Young's modulus $(\rho)$, the coefficient of thermal expansion $(\alpha)$ and the thermal conductivity $(\kappa)$ are: $E_c=$ 151e$^9$ N/m$^2$, $\alpha_c=$ 10e$^{-6}$/$^\circ$C, $\kappa_c=$ 2.09 W/mK for ZrO$_2$ and $E_m=$ 70e$^9$ N/m$^2$, $\alpha_m =$ 23e$^{-6}$/$^\circ$C, $\kappa_m =$ 204 W/mK for Al. Poisson's ratio $\nu$ is assumed to be constant and taken as 0.3 for the current study~\cite{zhaolee2009a}.  Here, the modified shear correction factor obtained based on energy equivalence principle as outlined in~\cite{Singh2011} is used. 

\paragraph*{Mechanical Buckling}
For mechanical buckling, we consider both uni- and bi-axial mechanical loads on the FGM plates. In all cases, we present the critical buckling parameters as, unless otherwise specified:

\begin{eqnarray}
\lambda_{\rm{uni}} = \frac{N_{\rm{xxcr}}^o b^2}{\pi^2 D_c} \nonumber \\
\lambda_{\rm{bi}} = \frac{N_{\rm{yycr}}^o b^2}{\pi^2 D_c} 
\label{eqn:mbuckparm}
\end{eqnarray}

where, $\lambda_{\rm{cru}}$ and $\lambda_{\rm{crb}}$ are the critical buckling parameters for uni- and bi-axial load, respectively,  $D_c=E_ch^3/(12(1-\nu^2))$. In order to be consistent with the literature, properties of the ceramic phase are used for normalization. 

\paragraph*{Thermal Buckling}
For thermal buckling, a temperature rise of $T_m=$ 5$^\circ$C in the metal-rich surface of the plate is assumed in the present study. In addition to nonlinear temperature distribution through the plate thickness, a linear distribution of the temperature is also considered in the present analysis by truncating the higher order terms in \Eref{eqn:heatconducres}. The plate is of uniform thickness and simply supported on all four edges.

\begin{rmk} In both the cases, the effect of the crack location and its orientation, the cutout radius, the plate thickness and boundary conditions on the global response is numerically studied. 
\end{rmk}

\paragraph{Validation} Table \ref{table:fgmHoleMBuck} shows the convergence of the critical buckling load with mesh size. Based on a progressive refinement, a 40$\times$40 structured quadrilateral mesh is found to be adequate to model the full plate for the present analysis. Before proceeding with a detailed parametric study on the effect of different parameters, the formulation developed herein is validated against available results pertaining to the critical buckling load of FGM plates with and without local defects. The computed critical buckling parameters for FGM plates with and without local defects under uniaxial mechanical loading are given in Tables \ref{table:isoMbuck} and \ref{table:fgmHoleMBuck}, respectively. It can be seen that the numerical results from the present formulation are in very good agreement with the existing solutions. 


\begin{table}[htpb]
\renewcommand\arraystretch{1.2}
\caption{Convergence study of critical buckling load $\lambda_{uni} = \frac{N_{xxcr}^o b^2}{\pi^2 D}$ for simply supported square Al/ZrO$_2$ plate with circular cut out at the center with $a/h=$ 100, $r/a=$ 0.1 for different gradient index.}
\centering
\begin{tabular}{lrrr}
\hline 
Mesh & \multicolumn{3}{c}{gradient index, $n$} \\
\cline{2-4}
& 0 & 1 & 5 \\
\hline
10 $\times$ 10 & 7.6439 & 5.3433 & 4.5458 \\
20 $\times$ 20 & 7.1399 & 4.9910 & 4.2461 \\
30 $\times$ 30 & 7.0640 & 4.9380 & 4.2009 \\
40 $\times$ 40 & 7.0247 & 4.9105 & 4.1760 \\
Ref.~\cite{zhaolee2009} & 6.9711 & 4.6858 & 4.0609 \\
\hline
\end{tabular}
\label{table:fgmHoleMBuck}
\end{table}

\begin{table}[htpb]
\renewcommand\arraystretch{1.2}
\caption{Comparison of critical buckling load $\lambda_{uni} = \frac{N_{xxcr}^o b^2}{\pi^2 D}$ for simply supported isotropic square plate subjected to uni-axial loading with $a/h=$ 10.}
\centering
\begin{tabular}{llrrr}
\hline
Boundary & & \multicolumn{3}{c}{$a/b$} \\
\cline{3-5}
&   & 1 & 1.5 & 2 \\
\hline
\multirow{2}{*}{SSSS} & Present & 3.7346 & 3.9814 & 3.7546 \\
      & Ref.~\cite{zhaolee2009a} & 3.7412 & 3.9613 & - \\
\cline{2-5}
\multirow{2}{*}{CCCC} & Present & 8.1818 & 6.8936 & 6.5347 \\
& Ref.~\cite{zhaolee2009a} & 8.1391 & 6.8891 & 6.5668 \\
\hline
\end{tabular}
\label{table:isoMbuck}
\end{table}

\subsection{Effect of crack} 
\begin{figure}[htpb]
\centering
\input{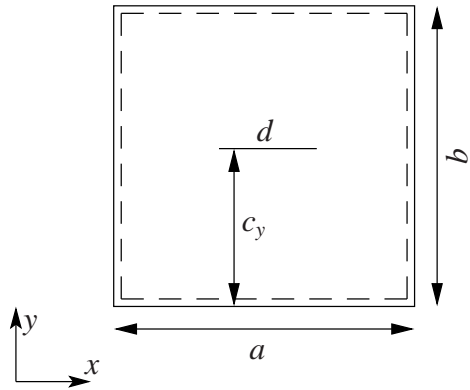}
\caption{Plate with a centrally located crack with simply supported boundary conditions.}
\label{fig:platedescribe}
\end{figure}

\begin{figure}[htpb]
\centering
\subfigure[Uni-axial]{\includegraphics[scale=0.75]{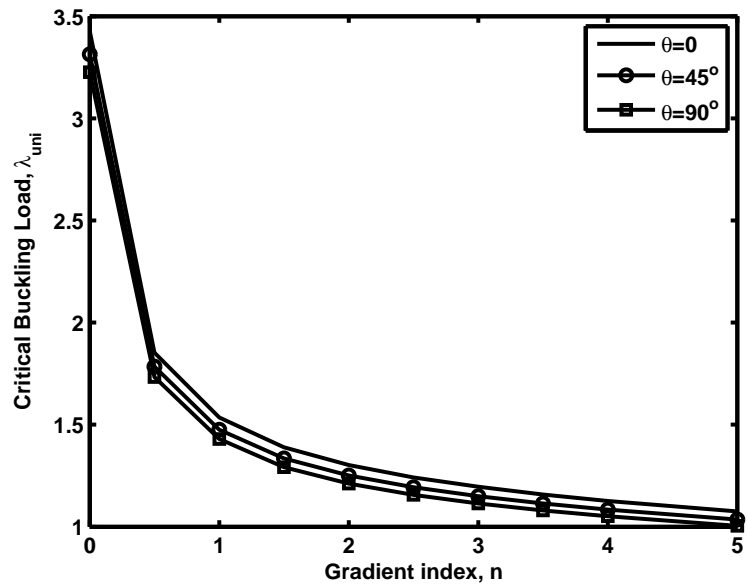}}
\subfigure[Bi-axial]{\includegraphics[scale=0.75]{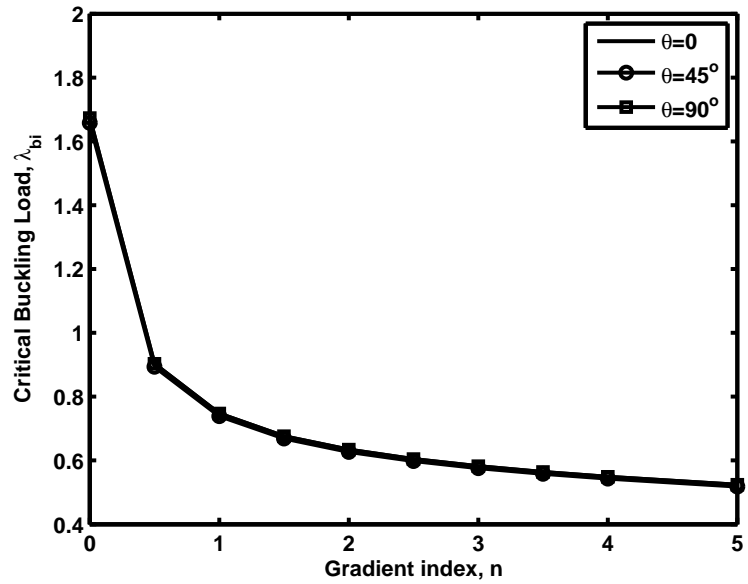}}
\caption{Critical buckling load as a function of gradient index for a simply supported square FGM plate subjected to uni-axial and bi-axial compressive loads with $a/h=$100 and crack length $d/a$=0.4.}
\label{fig:mbuck_nVari_SS2}
\end{figure}

\begin{figure}[htpb]
\centering
\subfigure[Uni-axial]{\includegraphics[scale=0.75]{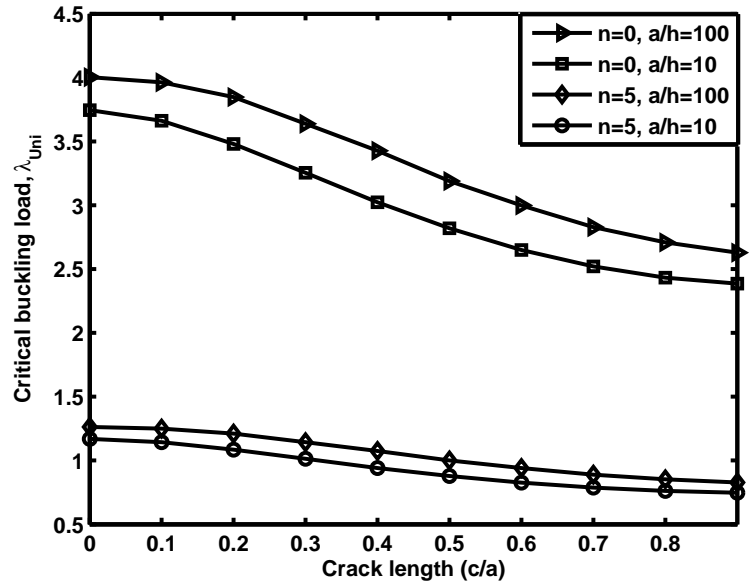}}
\subfigure[Bi-axial]{\includegraphics[scale=0.75]{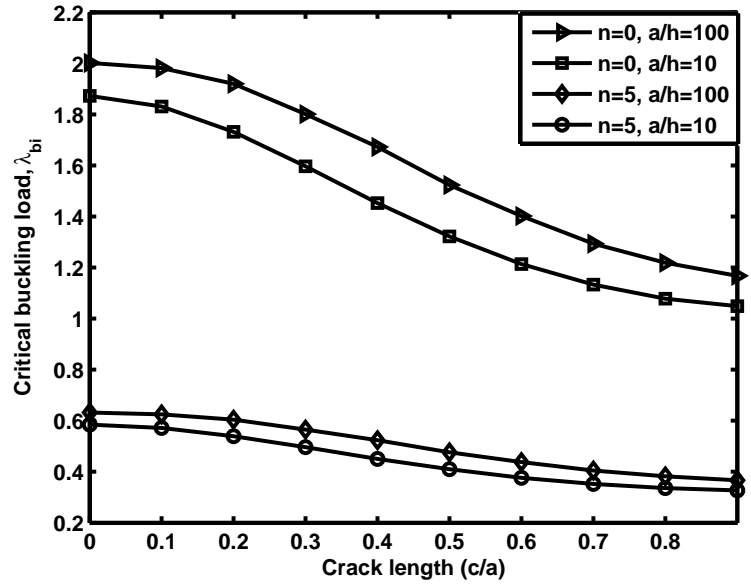}}
\caption{Critical buckling load as a function of crack length $c/a$ for a simply supported square FGM plate subjected to uni-axial and bi-axial compressive loads.}
\label{fig:mbuck_caVari_SS2}
\end{figure}

Consider a plate of uniform thickness, $h$ and with length and width as $a$ and $b$, respectively. \fref{fig:platedescribe} shows a with a center crack of length $d$ and at a distance of $c_y$ from the $x-$ axis with all edges simply supported. In this section, we study the influence of the crack length, the crack orientation, the gradient index and the type of loading on the critical buckling load for a square simply supported FGM plate with thickness $a/h=$ 10 and 100. \fref{fig:mbuck_nVari_SS2} shows the influence of the gradient index, $n$ and the crack orientation, $\theta$ on the critical buckling load for a FGM plate with a crack length $d/a=$ 0.4, subjected to uni- and bi-axial compressive loads. It is observed that with increase in gradient index $n$, the critical buckling load decreases for both the uni- and bi-axial compressive loads. This is because of the stiffness degradation due to increase in the metallic volume fraction. \fref{fig:mbuck_caVari_SS2} and Table \ref{tab:mechCrtLoad_Vs_angle_and_Gindex} shows the influence of the crack length and the crack orientation for a square FGM plate subjected to uni- and bi-axial compressive loads. It can be observed that increasing the crack length and the crack orientation decreases the critical buckling load. Further, it is observed that the frequency is lowest for a crack orientation $\theta=$ 90$^\circ$ when the plate is subjected to uni-axial compressive load. When the plate is subjected to bi-axial compressive load, the critical buckling load initially decreases, with further increase in crack orientation, the critical buckling load increases and reaches maximum at $\theta=$90$^\circ$. 

\begin{table}[htbp]
\centering
\renewcommand{\arraystretch}{1.5}
\caption{Influence of the crack orientation and the gradient index on the critical buckling load for a simply supported square FGM plates subjected to mechanical loads, viz., uni, biaxial compression with $a/h=$100 and crack length $c/a=$0.4.}
\begin{tabular}{lrrrrrrrr}
\hline
Crack & \multicolumn{8}{c}{gradient index, $n$} \\
\cline{2-9}
Angle,& \multicolumn{2}{c}{0} && \multicolumn{2}{c}{1} && \multicolumn{2}{c}{5}  \\
\cline{2-3}\cline{5-6}\cline{8-9}
$\beta$     & uni   & bi   & & uni   & bi    && uni   & bi \\
\hline
0     & 3.4268 & 1.6729 && 1.5348 & 0.7473 && 1.0749 & 0.5236 \\
30    & 3.3527 & 1.6549 && 1.4962 & 0.7378 && 1.0496 & 0.5177 \\
45    & 3.3141 & 1.6582 && 1.4753 & 0.7391 && 1.0352 & 0.5183 \\
90    & 3.2263 & 1.6729 && 1.4294 & 0.7473 && 1.0044 & 0.5236 \\
\hline
\end{tabular}%
\label{tab:mechCrtLoad_Vs_angle_and_Gindex}%
\end{table}%

\begin{figure}[htpb]
\centering
\subfigure[Linear temperature variation]{\includegraphics[scale=0.75]{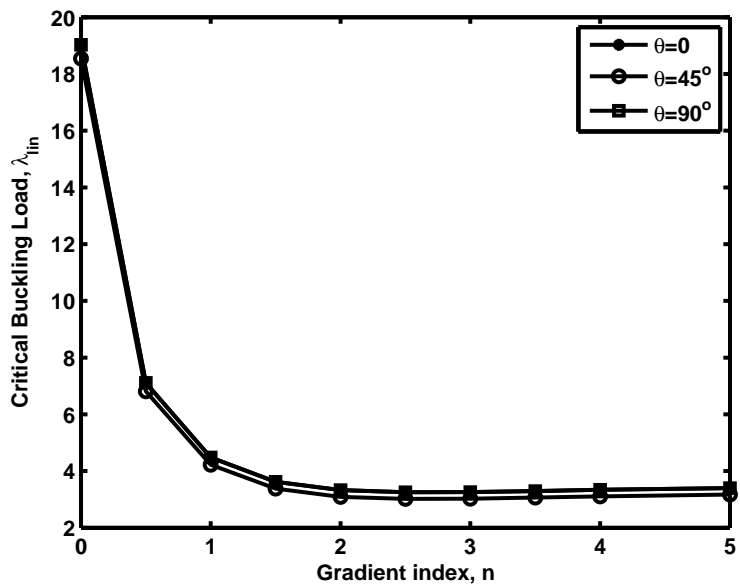}}
\subfigure[Non-linear temperature variation]{\includegraphics[scale=0.75]{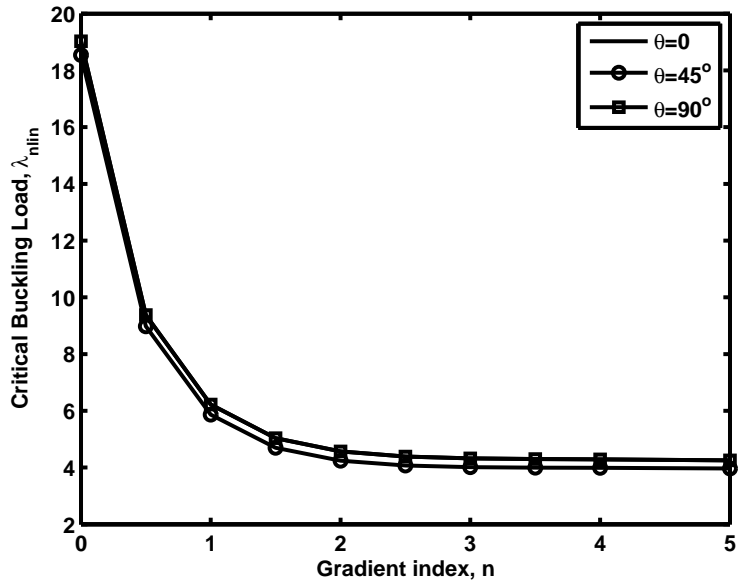}}
\caption{Critical buckling load as a function of gradient index for a simply supported square FGM plate subjected to linear and non-linear through thickness temperature variation with $a/h=$100 and crack length $c/a$=0.4.}
\label{fig:thbuck_nVari_SS2}
\end{figure}

The influence of the crack orientation and the gradient index on the critical buckling temperature for a FGM square plate with $a/h=$ 100, when subjected to linear and non-linear temperature rise through the thickness can be seen from Table \ref{tab:thCrtLoad_Vs_angle_and_Gindex} and also depicted in \fref{fig:thbuck_nVari_SS2}. It can be seen that increasing the gradient index decreases the critical buckling temperature both for the linear and non-linear temperature rise through the thickness. This behaviour can be attributed to the increasing metallic volume fraction when the gradient index increases. When the crack orientation is increased the critical buckling load decreases and reaches minimum when the crack orientation is in between $\beta=$ 30$^\circ$ and 45$^\circ$ and upon further increase, the critical buckling load increases. 

\begin{figure}[htpb]
\centering
\subfigure[Thick plate]{\includegraphics[scale=0.75]{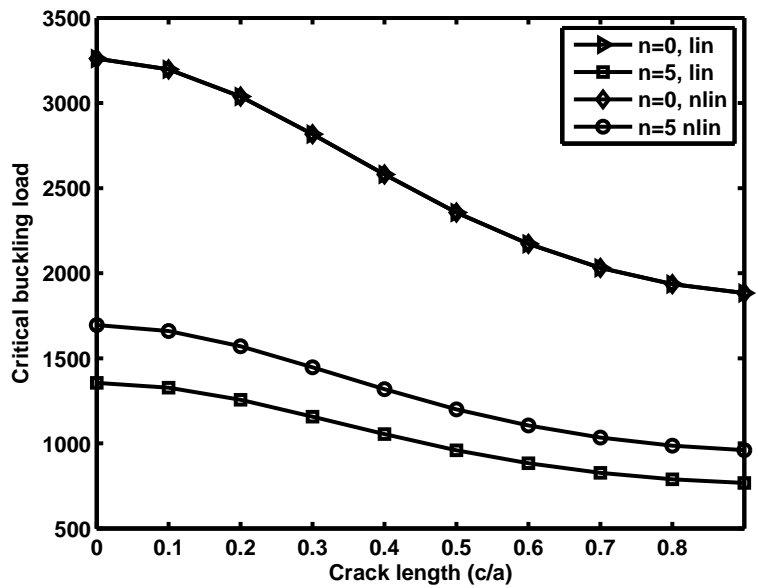}}
\subfigure[Thin Plate]{\includegraphics[scale=0.75]{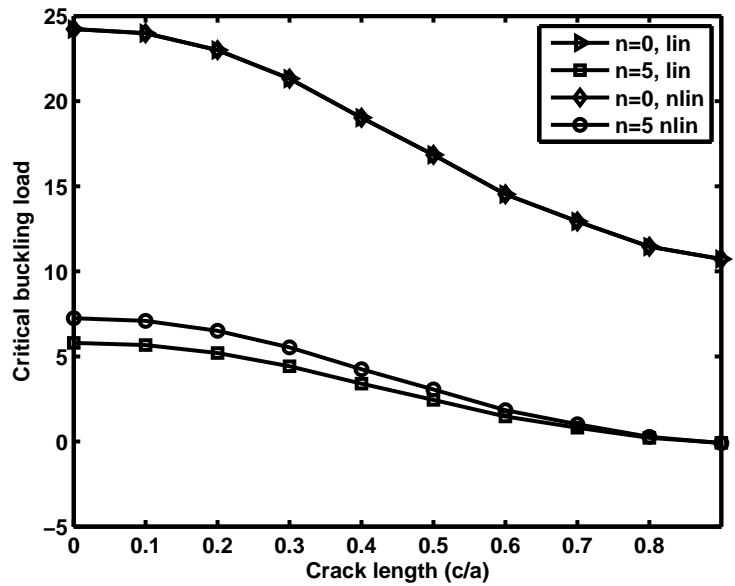}}
\caption{Critical buckling load as a function of crack length $c/a$ for a simply supported square FGM plate with thickness $a/h=$ 100 and 10 subjected to linear and non-linear through thickness temperature variation.}
\label{fig:thbuck_caVari_SS2}
\end{figure}

\begin{figure}[htpb]
\centering
\subfigure[Mechancial]{\includegraphics[scale=0.75]{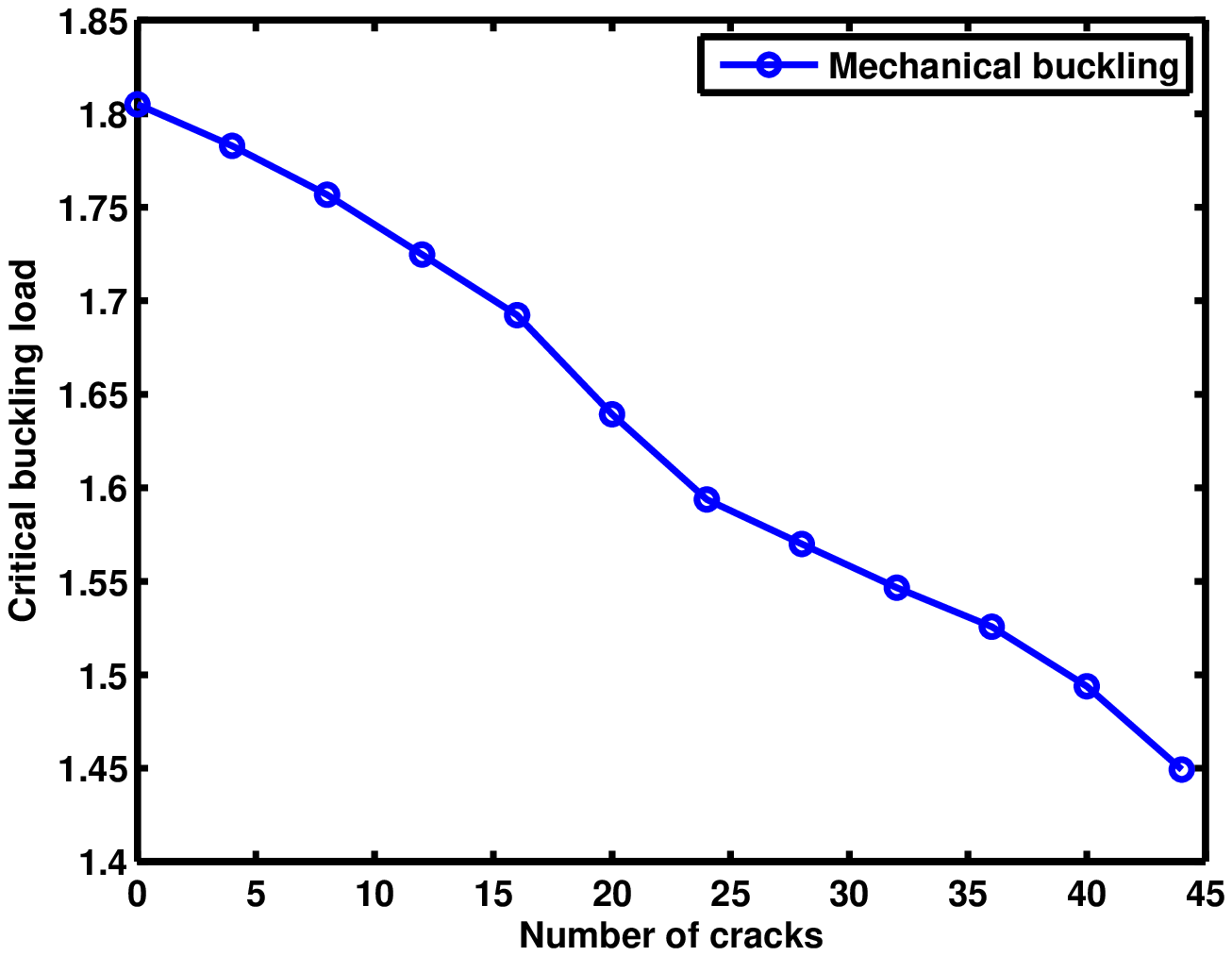}}
\subfigure[Thermal]{\includegraphics[scale=0.75]{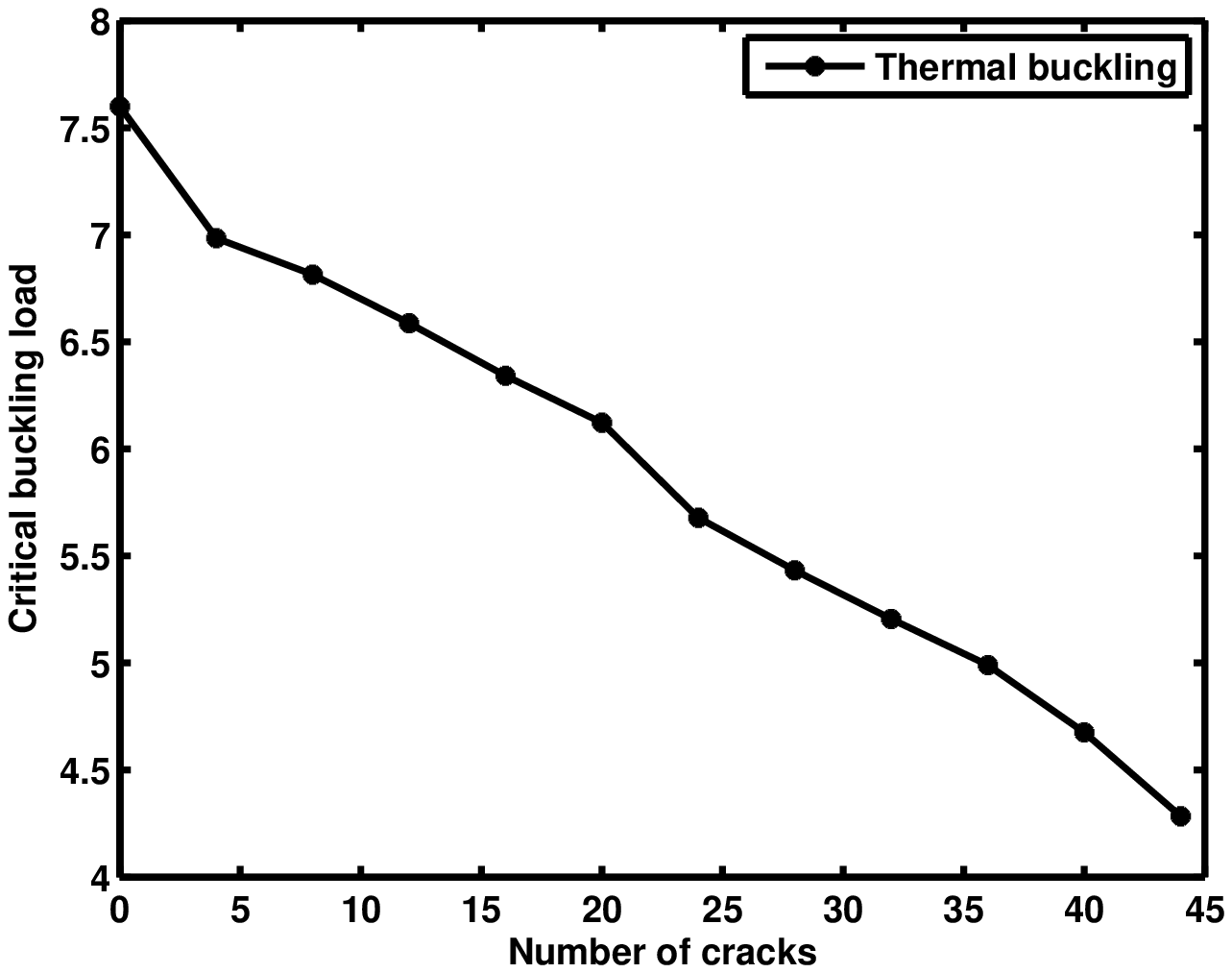}}
\caption{Critical buckling load as a function of number of cracks for a simply supported square FGM plate with $a/h=$100 and $n=$1 subjected to (a) uni-axial compressive loading (b) linear through thickness temperature variation.}
\label{fig:multi_crack}
\end{figure}

The influence of the crack length $d/a$ and the gradient index $n$ on the critical buckling temperature for a simply supported FGM square plate with $a/h=$ 10 and 100, when subjected to linear and non-linear temperature rise through the thickness is shown in \fref{fig:thbuck_caVari_SS2}. It can be seen that the combined effect of increasing the crack length, the gradient index and the plate thickness is to decrease the critical buckling temperature for both linear and non-linear temperature rise through the thickness. The combined effect can be attributed to the stiffness degradation due to the increase in the metallic volume fraction and the geometry change. The influence of temperature distribution through the thickness can be clearly seen from \fref{fig:thbuck_nVari_SS2}. It should be noted that there is no effect of the temperature distribution on the critical buckling temperature when the gradient index $n=$ 0 as expected, because the gradient index $n=$ 0, corresponds to pure ceramic plate. This holds true for both thick and thin plates. \fref{fig:multi_crack} shows the influence of number of cracks on the critical buckling load and critical buckling temperature for a square simply supported FGM plate with $a/h=$ 100 and gradinet index $n=$ 1. It can be seen that with increasing number of cracks, the critical buckling load and the critical buckling temperature decreases due to the increase in the local flexibility because of the presence of the discontinuity surface.
\begin{table}[htbp]
\renewcommand{\arraystretch}{1.5}
\centering
\caption{Influence of the crack orientation and the gradient index on the critical buckling load for a simply supported square FGM plates subjected to thermal loads, viz., through thickness temperature variation with $a/h=$100 and crack length $c/a=$0.4.}
\begin{tabular}{rrrrrrrrr}
\hline
Crack & \multicolumn{8}{c}{gradient index, $n$} \\
\cline{2-9}
Angle, & \multicolumn{2}{c}{0} & & \multicolumn{2}{c}{1} & & \multicolumn{2}{c}{5}  \\
\cline{2-3}\cline{5-6}\cline{8-9}
$\beta$ & T$_{_L}$   & T$_{_{NL}}$   & & T$_{_L}$   & T$_{_{NL}}$   & & T$_{_L}$   & T$_{_{NL}}$ \\
\hline
0     & 19.0289 & 19.0289 && 3.2929 & 4.5708 && 3.3978 & 4.2494 \\
30    & 18.4965 & 18.4965 && 3.0552 & 4.2410 && 3.1570 & 3.9483 \\
45    & 18.5432 & 18.5432 && 3.0791 & 4.2741 && 3.1689 & 3.9631 \\
90    & 19.0289 & 19.0289 && 3.2929 & 4.5708 && 3.3978 & 4.2494 \\
\hline
\end{tabular}%
\label{tab:thCrtLoad_Vs_angle_and_Gindex}%
\end{table}%

\subsection{Effect of cutouts} 
\begin{figure}[htpb]
\centering
\includegraphics[scale=0.7]{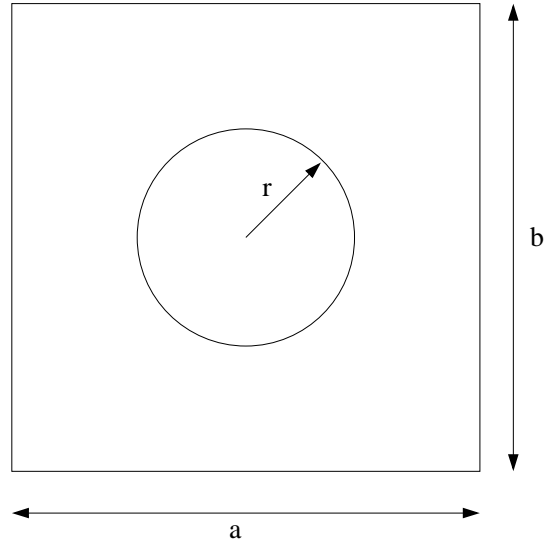}
\caption{Plate with a centrally located circular cutout. $r$ is the radius of the circular cutout.}
\label{fig:phole}
\end{figure}

\begin{figure}[htpb]
\centering
\includegraphics[scale=0.75]{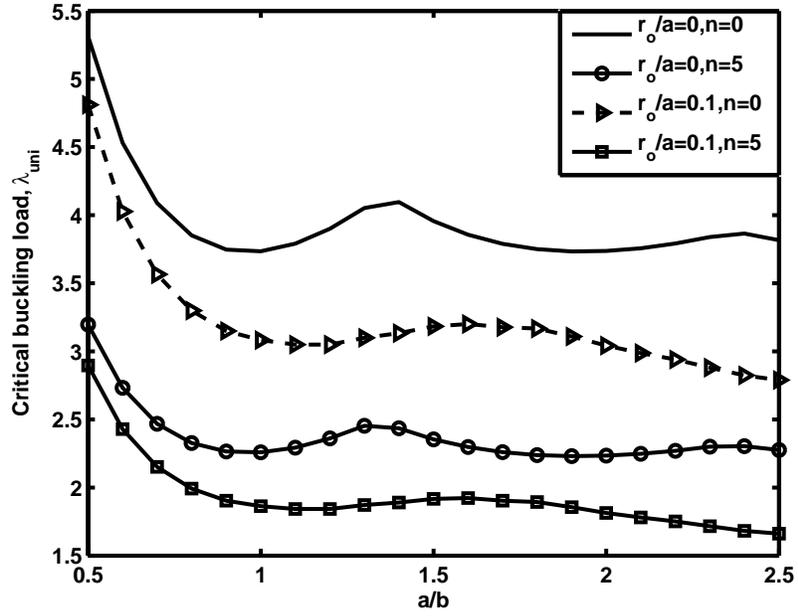}
\caption{Influence of the plate aspect ratio $a/b$ on the critical buckling load, $\Lambda_{uni} = \frac{N_{xxcr}^o b^2}{\pi^2 D}$, of the FGM plate subjected to uniaxial loading for different gradient index $n$ with a center cutout of radius, $r_o/a=$ 0.1 and $h/b=$ 0.1.}
\label{fig:MechBuck_ab_roaVari}
\end{figure}

\begin{figure}[htpb]
\centering
\includegraphics[scale=0.75]{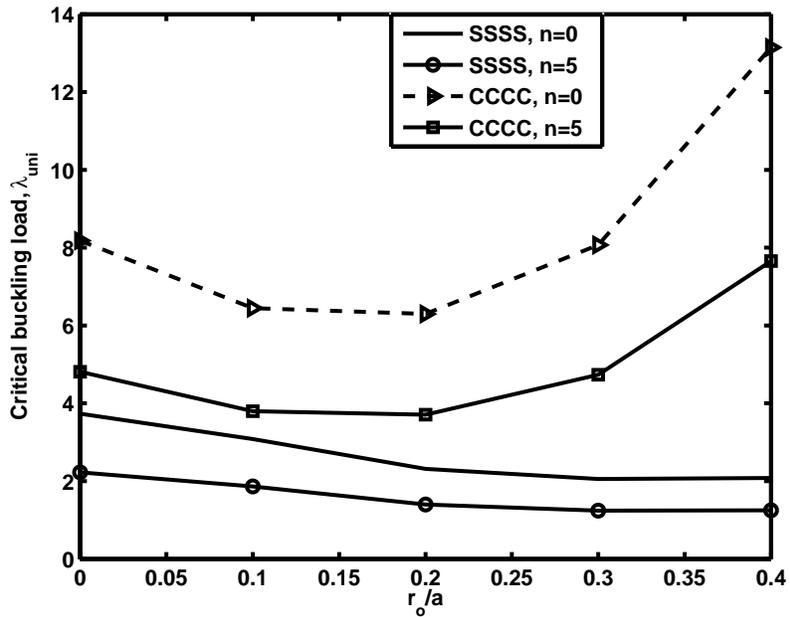}
\caption{Variation of the critical buckling load, $\Lambda_{uni} = \frac{N_{xxcr}^o b^2}{\pi^2 D}$ with cutout dimensions for a square FGM plate with $a/h=10$ subjected to uniaxial loading for different gradient index $n$ and boundary conditions.}
\label{fig:MechBuck_holeRadVari}
\end{figure}

\begin{figure}[htpb]
\centering
\subfigure[Mechancial]{\includegraphics[scale=0.75]{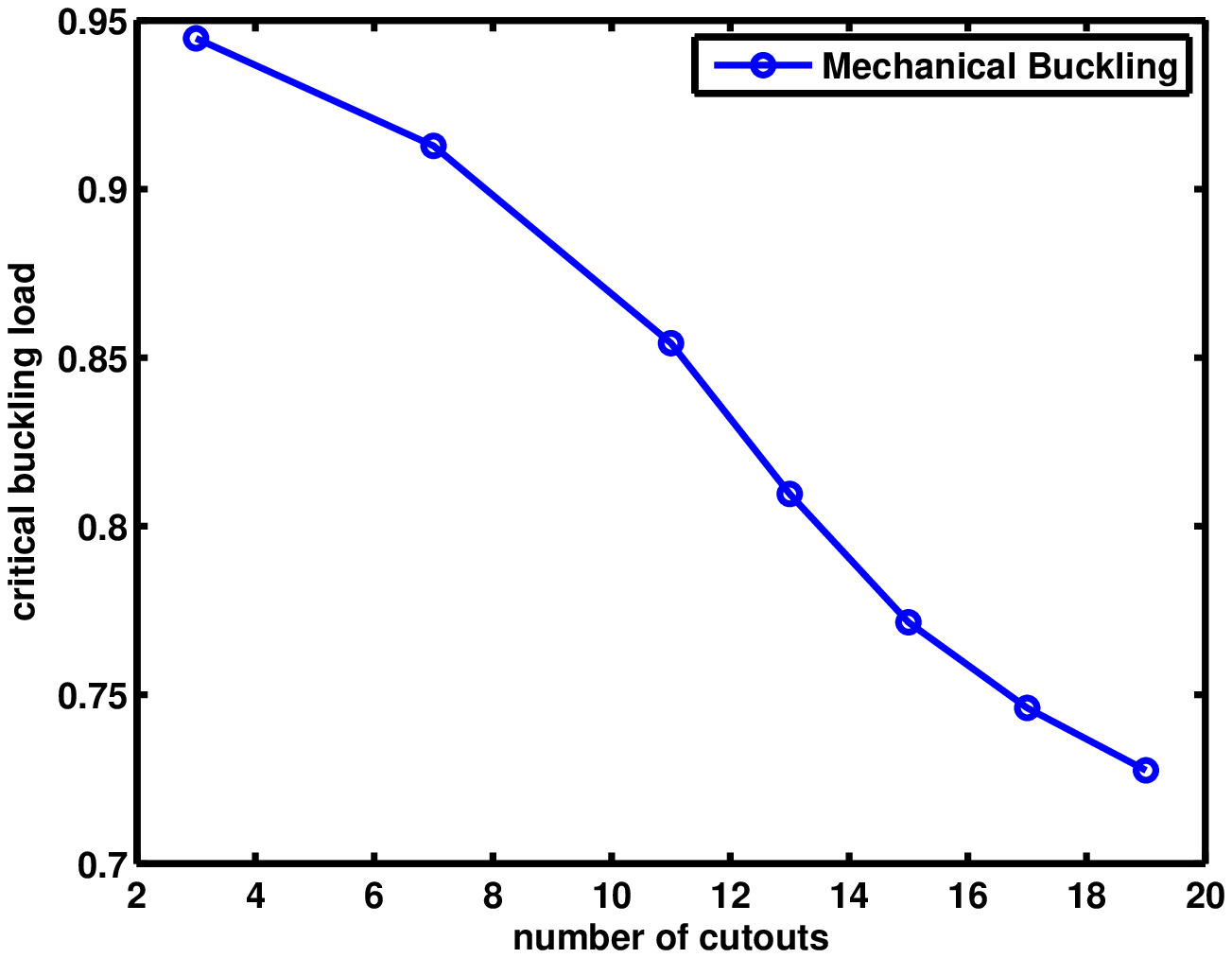}}
\subfigure[Thermal]{\includegraphics[scale=0.75]{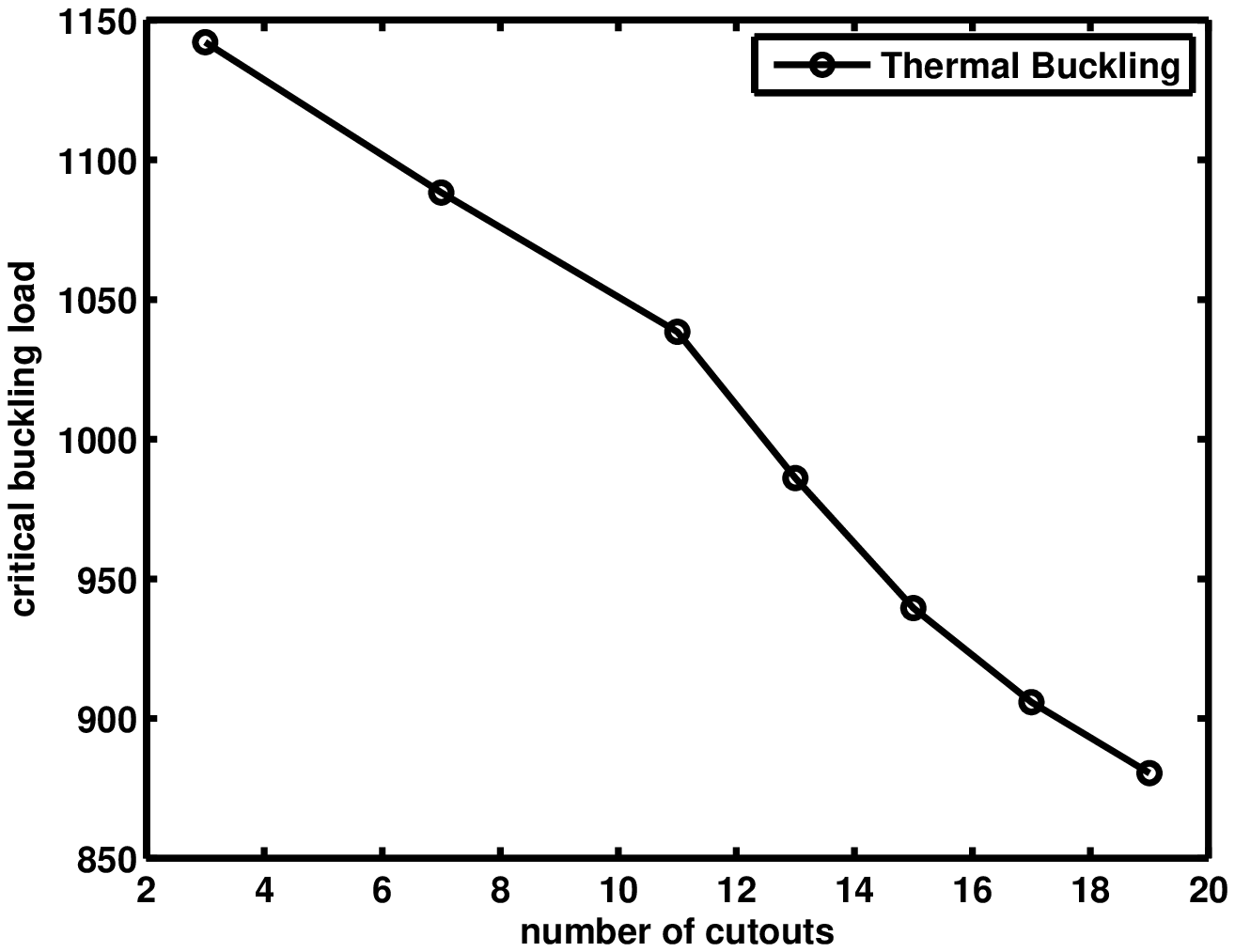}}
\caption{Critical buckling load as a function of number of cutouts for a simply supported square FGM plate with $a/h=$100 and $n=$5 subjected to (a) uni-axial compressive loading (b) nonlinear through thickness temperature variation.}
\label{fig:multi_cut_outs}
\end{figure}

Next, we study the influence of a circular cutout on the critical buckling load and critical buckling temperature for a FGM plate.  \fref{fig:phole} shows the geometry of the plate with a centrally located circular cutout. \fref{fig:MechBuck_ab_roaVari} shows the influence of the plate aspect ratio, the gradient index and a centrally located cutout on the critical buckling load when the FGM plate is subjected to a uniaxial compressive load. It can be seen that the combined effect of increasing the gradient index, the cutout radius and the plate aspect ratio is to decrease the critical buckling load. The decrease in the critical load due to the gradient index and the cutout radius can be attributed to the stiffness degradation due to increased metallic volume fraction and due to the presence of a discontinuity, respectively. The influence of the plate aspect ratio is attributed to the geometry effect. \fref{fig:MechBuck_holeRadVari} shows the influence of the centrally located circular cutout and the gradient index on the critical buckling load under two different types of boundary conditions, viz., all edges simply supported and all edges clamped. Here again, the plate is subjected to a uni-axial compressive load along the $y-$ direction. It can be seen that increasing the gradient index decreases the critical buckling load due to increasing metallic volume fraction, whilst, increasing the cutout radius lowers the critical buckling load in the case of simply supported boundary conditions and increases in the case of all edges clamped boundary. The influence of the boundary condition on the critical buckling load can be clearly observed from this study. The influence of randomly located circular cutout on the critical buckling load and the critical buckling temperature is shown in \fref{fig:multi_cut_outs}. A simply supported square FGM plate with $a/h=$ 100 and gradient index, $n=$ 5 is used and is subjected to uni-axial compressive load and nonlinear temperature rise through the thickness of the plate. It can be seen that increasing the number of cutouts, decreases the critical buckling load and the critical buckling temperature. This can be attributed to stiffness degradation due to the presence of cutouts.

\section{Conclusion}
In this paper, we did a systematic parametric study to bring out the influence of local defects, viz., cracks and cutouts on FGM plates when subjected to in-plane compressive load or linear/non-linear temperature distribution through the thickness of the plate. The plate kinematics is based on first order shear deformation theory. From the detailed numerical study, it is concluded that:

\begin{itemize}
\item Among all  other parameters the gradient index has more profound influence on critical buckling load of the FGM plate irrespective of plate geometry, support condition, thickness and crack/cutout geometry. Increasing the material gradient index $n$, decreases the buckling load. This is due to the increase in the metallic volume fraction.
\item Increasing the crack length or the cutout radius decreases the critical buckling load of the FGM plate.
\item Increasing the number of cracks/cutouts, decreases the overall stiffness of the plate and thus decreases the critical buckling load.
\end{itemize}

\section*{Acknowledgement}
S Natarajan would like to acknowledge the financial support of the School of Civil and Environmental Engineering, The University of New South Wales for his research fellowship since September 2012.

\bibliography{FGM_buckling}
\bibliographystyle{elsarticle-num}

\end{document}